\documentstyle[12pt]{article}
\newcommand{\vs}{\vspace{.2in}}

\def \ind{{\rm ind\,}}
\def \tip{{\rm tip\,}}

\def \p{\partial}
\def\C{{\bf C}}

\def\Z{{\bf Z}}
\def \Q{{\bf Q}}
\def \P{{\bf P}}
\def\G{{\Gamma}}

\def \s{\sigma}

\def \G{\Gamma}
\def \~ {\tilde}
\def \D{{\~ D}}
\def \_{\bar}
\def \sp{\,\,\,\,\,\,\,}
\def \[{\lceil}
\def \]{\rceil}

\def \be{\begin{equation}}
\def \ee{\end{equation}}
\def \noin{\noindent}
\def\qed{\hfill $\Box$}

\begin{document}

\title{ On the number of singular points of plane curves
\footnote{ Mathematics
Subject Classification: 14H20, 14H10, 14D15, 14N05}}

\bigskip

\author{ S. Orevkov and  M. Zaidenberg }

\date{}
\maketitle

\bigskip
\begin{abstract} This is an extended, renovated and updated report on our
joint work [OZ]. The main result is an inequality for the numerical type of
singularities of
a plane curve, which involves
the degree of the curve, the multiplicities and the
Milnor numbers of its singular points. It is a corollary of the logarithmic
Bogomolov-Miyaoka-Yau's type  inequality due to Miyaoka.
It was first proven by F. Sakai at 1990
and rediscovered by the authors independently in the particular case of an
irreducible cuspidal curve at 1992. Our proof is based on the localization,
the local Zariski--Fujita decomposition and uses a graph discriminant calculus.
The key point is a local analog of the BMY-inequality for a plane curve germ.
As a corollary, a boundedness criterium for a family of plane curves has been
obtained. Another application of our methods is the following fact: a rigid
rational cuspidal plane curve cannot have more than 9 cusps.
\end{abstract}

\noin This is an extended, renovated and updated report on our joint work [OZ]
which
 the second named author presented at the Conference on Algebraic Geometry
 held at
 Saitama University, 15-17 of March, 1995. It is his pleasure to thank
 Professor F. Sakai who invited him to take part in this nice conference. \\

\section{Asymptotics of the number
of ordinary cusps}

We start with a brief survey of known results in the simplest case of ordinary
 cusps.

It is well known that for a nodal plane curve $D \subset \P^2$ of degree $d$
 the number of nodes can be an arbitrary non--negative integer allowed by the
 genus formula, i.e. any integer from the interval $[0,\,{d-1 \choose 2}]$.
If $D$ is a Pl\"ucker curve with only ordinary cusps as singularities, which
 has $\kappa$ cusps, then still $$\kappa \le {d-1 \choose 2}\,,$$ but this
 time the inequality is strict starting with $d = 5$. Indeed, by Pl\"ucker
 formulas $$0 < d^* = d(d-1) - 3\kappa$$ where $d^*$ is the class of $D$, and
 $$0 \le f = 3d(d-2) - 8\kappa$$ where $f$ is the number of inflexion points
 of $D$. Thus, we have
\be
\kappa < {1 \over 3}d(d-1)
\ee
and
\be \kappa \le {3 \over 8}d(d-2)\,,
\ee
which is strictly less than ${d-1 \choose 2}$ for $d \ge 5$.

 Note that $d \le 4$ for a rational cuspidal Pl\"ucker curve $D$, due to (2)
 and the genus formula. Therefore, up to projective equivalence there exists
 only two such curves, namely the cuspidal cubic and the Steiner
 three-cuspidal quartic (we suppose here that at least one cusp really occurs,
 otherwise we have to add also the line and the smooth conic).

{}From (2) it follows that
\be \limsup\limits_{d \to \infty} {\kappa \over d^2} \le {3 \over 8} \,\,.\ee
Using the spectrum of singularity (or, equivalently, the Mixed Hodge
 Structures) A. Varchenko [Va] found an estimate
\be \limsup\limits_{d \to \infty} {\kappa \over d^2} \le {23 \over 72} \,\,,\ee
which is better than (3) by ${1 \over 18}$.
Another ${1 \over 144}$ was gained in the work of F. Hirzebruch and T.
 Ivinskis [H, Iv] by applying Miyaoka's logarithmic form of the
 Bogomolov-Miyaoka-Yau (BMY) inequality:
\be \limsup\limits_{d \to \infty} {\kappa \over d^2} \le {5 \over 16} \,\,.\ee
Furthermore, in this work an elegant example was given which shows that
\be \limsup\limits_{d \to \infty} {\kappa \over d^2} \ge {1 \over 4} \,,\ee
where $c$ means now the maximal number of cusps among all the cuspidal
 Pl\"ucker curves of degree $d$. \\

\noin {\bf Example} [H, Iv]. Starting with a generic smooth cubic $C$,
 consider its dual curve $D = C^*$, which is an elliptic sextic with nine
 ordinary cusps as the only singularities. Let $F = 0$ be the defining
 equation of $D$. Set $D_k = \{F(x^k : y^k : z^k) = 0\}$. Then $D_k$ is again
 a cuspidal Pl\"ucker curve. It has degree $d_k = 6k$ and $\kappa_k = 9k^2$
 cusps (indeed, $(x : y : z) \longmapsto  (x^k : y^k : z^k)$ is a branched
 covering $\P^2 \to \P^2$ of degree $k^2$ ramified along the coordinate axes
 which meet $D$ normally). Thus, here $\kappa_k = d_k^2 / 4$.

In fact, the lower bound $1/4$ can be improved, by a similar method, at least
 by $1/32$ (A. Hirano [Hi]). Together with (5) this yields
 \be {10 \over 32} \ge \limsup\limits_{d \to \infty} {\kappa \over d^2} \ge {9
 \over 32}\,, \ee
which is still far away from giving the exact asymptotic. See also [Sa] for a
 discussion on what is known for small values of $d$.

\section{The main inequality. Bounded families of plane curves}

Next we consider, more generally, plane curves with arbitrary singularities.
 By {\it a cusp} we mean below a locally irreducible singular point. We say
 that $D \subset \P^2$ is {\it a cuspidal curve  } if all its singular points
 are cusps. The following theorem, which is the main result presented in the
 talk, was first proven by F. Sakai [Sa]. Independently and later it was also
 found by the authors in the special case of cuspidal curves [OZ] (actually,
 the proof in [OZ] goes through without changes for nodal--cuspidal curves,
 i.e. plane curves with nodes as the only reducible singular points). Both
 proofs are based on the logarithmic version of the BMY-inequality due to
 Miyaoka [Miy], but technically they are different. \\

\noin {\bf Theorem 1.} {\it Let $D$ be a plane curve of degree $d$ with the
 singular points $P_1, \dots, P_s$. Let $\mu_i$ resp. $m_i$ be the Milnor
 number resp. the multiplicity of $P_i \in D$. If $\P^2 \setminus D$ has a
 non--negative logarithmic Kodaira dimension, then
\be \sum_{i=1}^s (1 + {1\over 2m_i})\,\mu_i \le d^2 - {3\over2}d\,\,. \ee
In particular,
 \be    \sum_{i=1}^s \mu_i \le {2m\over 2m+1}(d^2 - {3\over 2}d) \,\, , \ee
where $m = \max\limits_{1\le i \le s} \{m_i\}$.} \\

\noin {\bf Remarks. } a) For $ m \le 3$ and $D$ irreducible Theorem 1 had
 been proved by Yoshihara [Yo1,2], whose work stimulated the later progress. \\

\noin b) For an irreducible plane curve $D$ of degree $d \ge 4$ the
 logarithmic Kodaira dimension $\bar\kappa(\P^2- D)$ is non--negative besides
 the case when $D$ is a rational cuspidal curve with one cusp; see [Wak]. \\

\noin {\bf Corollary 1. } {\it If $D \subset \P^2$ is an irreducible cuspidal
 curve of geometric genus $g$, then under the assumptions of Theorem 1 one has
 \be  g\ge {d^2-3(m+1)d\over 2(2m+1)}+1 \,\, . \ee
In particular, a family of such curves is bounded iff $g$ and $m$ are bounded
 throughout the family. } \\

In general, the latter conclusion does not hold for non-cuspidal curves
 (indeed, the family of all the irreducible rational nodal plane curves is
 unbounded). However, it becomes true if one replaces the geometric genus $g =
 g(D)$ by the Euler characteristic $e(D)$ (thus, involving not only the
 topology of the normalization, but the topology of the plane curve itself).
 Moreover, in this form it works even for reducible curves. \\

\noin {\bf Corollary 2. } {\it Under the assumptions of Theorem 1 one has
\be d(d-3(m+1)) \le (2m+1)(-e(D)) \,\,.\ee
Therefore, a family of (reduced) plane curves is bounded iff the absolute
 value of the Euler characteristic and the maximal multiplicity of the
 singular points are bounded throughout the family. } \\

The Corollary easily follows from Theorem 1 and the formula (see [BK])
 $$\sum_{i=1}^s \mu_i = d(d-3) + e(D)\,\,.$$

In the case of irreducible curves, in the estimate (11) it is convinient to
 use the first Betti number $b_1(D) = 2 - e(D)$ instead of the Euler
 characteristic. In particular, for irreducible nodal curves it is the only
 parameter involved.

An immediate consequence of (11) is that $d \le 3m+3$ if $e(D) \ge 0$.
 Furthermore, $d \le 3m+2$ if $e(D) > 0$; this is so, for instance, if $D$ is
 a rational cuspidal curve. In fact, in the latter case $d < 3m$ [MaSa], and
 also by the genus formula $$\sum\limits_{i=1}^s {\mu_i \over m_i} \le 3d-4$$
 [OZ].

\section{BMY-inequalities}

These inequalities provide the basic tool in the proof of Theorem 1. Let
 $\s \,:\,X \to \P^2$ be the minimal embedded resolution  of singularities of
 $D$, and let ${\~ D} \subset X$ be the reduced total preimage of $D$. Thus,
 ${\~ D}$ is a reduced divisor of simple normal crossing type, and
 ${\~ D}=\s^{-1}(D)$. Let $K = K_X$ be the canonical divisor of $X$. If
 ${\bar k}(\P^2 \setminus D) = {\bar k}(X \setminus {\~ D}) \ge 0$, i.e. if
 $|m(K+{\~ D})| \neq 0$ for $m$ sufficiently large, then (see [Fu]) there
 exists {\it the Zariski decomposition}  $K +  \D = H + N$, where $H,\,N$ are
 $\Q$--divisors in $X$ such that \\

\noin i) the intersection form of $X$ is negatively definite on the subspace
 $V_N \subset {\rm Pic}X \otimes \Q$ generated by the irreducible componenets
 of $N$; \\

\noin ii) $H$ is nef, i.e. $HC \ge 0$ for any complete irreducible curve
 $C \subset X$;\\

\noin iii) $H$ is orthogonal to the subspace $V_N$.\\

\noin By (iii) we have $$(K+{\~ D})^2 = H^2 + N^2\,,$$ where $N^2 \le 0$.
 Thus, $H^2 \ge  (K+{\~ D})^2$. \\

\noin {\bf Theorem} (Y. Miyaoka [Miy];
  R. Kobayashi--S.
  Nakamura--F. Sakai [KoNaSa]). \\
{\it
\noin a) If ${\bar k}(\P^2 \setminus D) \ge 0$, then
\be (K+{\~ D})^2 \le 3e(\P^2 \setminus D)\,.\ee
b) If ${\bar k}(\P^2 \setminus D) = 2$, then
\be H^2 \le 3e(\P^2 \setminus D)\,.\ee}

\noin {\bf Remark.} (13) holds, for instance, in the case when $D$ is an
 irreducible curve  with at least three cusps [Wak]. \\

Next we describe an approach to the proof of Theorem 1, mainly following [OZ].
  An advantage of this approach is that, in the particular case of irreducible
 cuspidal curves, we obtain formulas which express all the ingradients of the
 above BMY-inequalities in terms of the Puiseux characteristic sequences of
 the cusps. In fact, we prove a local version of Theorem 1 for the case of a
 cusp (see Theorem 2 in Section 11 below). Together with the  BMY-inequality
 (12) this provides a proof of Theorem 1 in the cuspidal case. A similar local
 estimate participates in the proof in [Sa], which is, by the way, much
 shorter. Instead of the Puiseux data it deals with the multiplicity sequences
 of the singular points. Combining both approches, we give in Section 12 a
 proof of the local estimate for arbitrary singularity (Theorem 3), thus
 proving Theorem 1 in general case. Actually, this is the proof of F. Sakai,
 with more emphasize separately on the local and the global aspects.

In the final Section 13 we apply the methods developed in the previous
 sections for pushing forward in the rigidity problem for rational cuspidal
 plane curves (see [FZ1,2]).

\section{Localization}

Let, as above, $D \subset \P^2$ be a plane curve of degree $d$ and let
 $\s\,:\,X \to \P^2$ be the  minimal embedded resolution of the singular
 points $P_1,\dots,P_s$  of $D$. Let $D'$ be the proper preimage of $D$ in
 $X$, ${\~ D} = D' \cup E$ be the reduced total preimage of $D$ and
 $E = E_1 \cup \dots \cup E_k, \,\,E_i = \s^{-1}(P_i),$ be the exceptional
 divisor of $\s$. Let $E_i = \sum\limits_{j=1}^{k_i} E_{ij}$ be the
 decomposition of $E_i$ into irreducible components. Fix also a line
 $L \subset \P^2$ which meets $D$ normally; denote by ${L'}$ the proper
 preimage of $L$ in $X$. Then, clearly, $\{E_{ij}\}$ and $L'$ form a basis of
 the vector space ${\rm Pic}X \otimes \Q$. Let $V_i = V_{E_i}$ be the subspace
 generated by the irreducible components $E_{ij}$ of $E_i$ and $V_{L'}$ be the
 one--dimensional subspace generated by $L'$ in ${\rm Pic}X \otimes \Q$. Since
 the intersection form of $X$ is non-degenerate, we have the orthogonal
 decomposition $${\rm Pic}X \otimes \Q = V_{L'} \oplus
 (\bigoplus_{i=1}^s V_i)\,.$$ Therefore, for each $i=1,\dots,s$ there exists
 the unique orthogonal projection ${\rm Pic}X \otimes \Q \to V_i$, and also
 such a projection onto the line $V_{L'}$. For any $\Q$--divisor $Z$ denote by
 $Z_{L'}$ resp. $Z_i$ its projection into $V_{L'}$ resp. into $V_i$. Then we
 have $$Z^2 = Z_{L'}^2 + \sum\limits_{i=1}^s Z_i^2\,.$$ In particular, since
 $K_{L'} + {\~ D}_{L'} = (d-3){L'}$ we have $$(K + {\~ D})^2 = (d-3)^2 +
 \sum\limits_{i=1}^s (K_i + {\~ D}_i)^2\,,$$ where the summands in the last
 sum are all negative (indeed, $E_i$ being an exceptional divisor, the
 intersection form of $X$ is negatively definite on the subspace $V_i$). It
 is easily seen that (8) follows from (12) and the local estimates
\be -(K_i + {\~ D}_i)^2 \le  (1 - {1 \over m_i}) \mu_i \,. \ee

In what follows we trace a way of proving (14). This is done in particular
 case of an irreducible singularity in Section 11 (Theorem 2), and in general
 in Section 12 (Theorem 3). Note that the assumption of local irredubicibility
 is important only in Sections 9, 10, 11.

\section{Weighted dual graph}

Let $E = E_1 \cup \dots \cup E_k$ be a curve with simple normal crossings in
 a smooth compact complex surface $X$. Assume, for simplicity, that all the
 irreducible components $E_i$ of $E$ are rational curves and that their
 classes in ${\rm Pic}X \otimes \Q$ are linearly independent. Let $A_E$ be the
 matrix of the intersection form of $X$ on the subspace $V_E =
 {\rm span}\,(E_1,\dots, E_k) \subset {\rm Pic}X \otimes \Q$ in the natural
 basis $E_1,\dots, E_k$ (we denote by the same letter a curve and its class in
 ${\rm Pic}X \otimes \Q$). Then $A_E$ is at the same time the incidence matrix
 of {\it the dual graph} $\Gamma_E$ of $E$, which is defined as follows. The
 vertices of $\G_E$ correspond to the irreducible components $E_i$ of $E$; two
 vertices $E_i$ and $E_j$, where $i \neq j$, are joint by a link
 $[E_i, \,E_j]$ iff $E_i\cdot E_j > 0$. The weight of the vertex $E_i$ is
 defined to be the self--intersection index $E_i^2$.

Let $C$ be another curve in $X$ which meets $E$ normally. Then we consider
 also {\it the dual graph $\G_{E,\,C}$ of $E$ near $C$}; it is the graph
 obtained from $\G_E$ by attaching $E_i\cdot C$ arrowheads to the vertex
 $E_i,\,\,i=1,\dots,k$. We denote by $\nu_i$ resp. ${\~ \nu}_i$ {\it the
 valency} of $E_i$ in $\G_E$ resp. in $\G_{E,\,C}$.

By {\it a twig} of a graph $\G$ one means an extremal linear branch of $\G$;
its end point is called {\it the tip} of the twig.

\section{Local Zariski--Fujita decomposition}

Since in the sequel we are working only locally over a fixed singular point
 $P=P_i$ of $D$, we change the notation. Omitting subindex $i$, from now on we
 denote by $E$ the corresponding exceptional divisor $E_i$ and by $V_E$ the
 corresponding subspace $V_i$. Thus, $K_E, \,{\~ D}_E, \,D'_E$ etc. mean the
 projections $K_i,\,{\~ D}_i,\,D'_i$... of the corresponding divisors into
 $V_E = V_i$. Set also $\mu = \mu_i$ and $m = m_i$. Note that in this case the
 dual graph $\G_E$ is a tree.

By {\it the local Zariski--Fujita decomposition} we mean the decomposition
 $$K_E + {\~ D}_E = H_E + N_E\,,$$ where $H_E,\,N_E \in V_E$ are effective
 $\Q$--divisors such that \\

\noin i) the support of $N_E$ coincides with the union of all the twigs of
 $\G_E$ which are not incident with the proper preimage $D'$ of $D$ in $X$,
 i.e. all the twigs of $\G_{E,\,D'}$ without arrowheads, and \\

\noin ii) $H_E$ is orthogonal to each irreducible component of
 ${\rm supp}\,N_E$. \\

Note that all the twigs in the ${\rm supp}\,N_E$ are {\it admissible}, i.e.
 all their weights are $\le -2$. Using non-degeneracy of the intersection
 form on an admissible twig, T. Fujita [Fu, (6.12)] proved that there exists
 the unique such decomposition. Moreover, he proved that up to certain
 exceptions the global Zariski decomposition $K + {\~ D} = H + N$ provides the
 local one via the projection (see [Fu, (6.20-6.24); OZ, Theorem 1.2]). Here
 we do not use this result, and so we do not give its precise formulation.

What we actually use is the equality $$(K_E + {\~ D}_E)^2 = H_E^2 + N_E^2\,.$$
 According to [Fu, (6.16); OZ, 1.1, 2.4], the latter summands can be computed
 in terms of the weighted graph $\G_{E,\,D'}$. This is done in the next
 section.

\section{Graph discriminants and inductances}

By definition, {\it the discriminant} $d(\G)$ of a weighted graph $\G$ is
$\det(-A)$, where $A$ is the incidence matrix of  $\G$
(or the intersection matrix of $E$, if $\G=\G_E$). It is easily seen that
 $d(\G_E) = 1$, because in our case $E$ is a contractible divisor.

{\it The inductance} of a twig $T$ of $\G$ is defined as
$$
        \ind(T) = {d(T-\tip(T)) \over d(T)}\,.
$$
Denote by $T_1,\dots,T_k$ the twigs of $\G_E$ which are not incident with $D'$
, i.e. the twigs of $\G_{E,\,D'}$. Then we have \\

\noin {\bf Lemma 1} [Fu, (6.16)]. {\it

\noin a)
\be -N_E^2 = \sum\limits_{i=1}^k \ind(T_i)\,.\ee
b) Let $v_T$ be the first vertex of $T = T_i\,(1\le i \le k)$, i.e. the vertex
 of $T$ opposite to the tip of $T$. Then the coefficient of $v_T$ in
 the decomposition of the divisor $N_E$ is equal to $1/d(T)$.} \\

 Since the graph $\G_E$ is a tree, for given vertices $E_i$ and $E_j$ (not
 necessary distinct) there is the unique shortest path in $\G_E$ which joins
 them. Denote by $\G_{ij}$ the weighted graph obtained from $\G_E$ by deleting
 of this path together with the vertices $E_i$ and $E_j$ themselves (and, of
 course, with all their incident links). So, in general the graph $\G_{ij}$ is
 disconnected.

Let $B_E = (b_{ij}) = A_E^{-1}$ be the inverse of the intersection matrix
 $A_E$. The following formula can be easily obtained by applying the Cramer
 rule. \\

\noin {\bf Lemma 2} [OZ, (2.1)].
\be b_{ij} = -d(\G_{ij})\,\,.\ee

Recall that $\bar\nu_i$ resp. $\nu_i$ denotes the valency of the vertex $E_i$
 of the graph $\G_{E,\,D'}$ resp. $\G_E$ and $\mu = \mu_i$ denotes the Milnor
 number of the singular point $P=P_i \in D$. We have\\

\noin {\bf Lemma 3} [OZ, (2.4), (2.7), (4.1)]. {\it In the notation as above
\be H_E^2 = \sum\limits_{\bar\nu_i>2,\,\bar\nu_j>2} b_{ij} c_i c_j \,,\ee
where $$c_i = (\bar\nu_i-2)-\sum {1\over d(T_j)}$$ and the last sum is taken
 over all the twigs $T_j$ which are incident with the vertex $E_i$;
\be (K_E+E)^2 = -2-\sum_{i=1}^n b_{ii}(\nu_i-2)\,,\ee and \footnote{Whereas
 (17) is valid for any rational SNC-tree $E$ with non-degenerate intersection
 form and admissible twigs on a smooth surface, (18) and
 (19) are true only when $E$ is the exceptional divisor of the resolution of
 singularity of a plane curve germ.}
\be \mu = 1-\sum_{i,j} b_{ij}(\bar\nu_i-2)(\bar\nu_j-\nu_j)\,. \ee} The proof
 of (17) is based on the Adjunction Formula and Lemmas 1,2; (18)
is proven by induction on the number of blow-ups; (19) follows from the
 adjunction formula and the formula
$$\,\,\,\,\,\,\,\,\,\,\,\,\,
\,\,\,\,\,\,\,\,\,\,\,\,\,\,\,\,\,\,\,
\,\,\,\,\,\,\,\,\,\,\,\,\,\,\,
\,\,\,\,\,\,\,\,\,\,\,\,\,\,\,\,\,
\mu = 1 - D'_E(K_E + {\~ D}_E)
 \,.\,\,\,\,\,\,\,\,\,\,\,\,\,\,
\,\,\,\,\,\,\,\,\,\,\,\,\,\,\,\,\,
\,\,\,\,\,\,\,\,\,\,\,\,\,\,\,\,
\,\,\,\,\,\,\,\,\,\,\,\,\,\,\,\,
\,(19')$$

\section{Calculus of graph discriminants}

To use the formulas from the preceding section we have to compute the entries
 $b_{ij}$ of the inverse $B_E = A_E^{-1}$, where $A_E$ is the intersection
 matrix of the exceptional divisor $E$; that means to compute corresponding
 graph discriminats (see Lemma 2). This section provides us with the necessary
 tools. They were developed in the work of Dr\"ucker-Goldschmidt [DG] (cf.
 also [Ra, Ne, Fu, (3.6)]) and afterwords interpreted by S. Orevkov [OZ] in
 the following way.

Let $\G$ be a weighted graph and let $A = A_{\G}$ be its intersection form.
 Recall that $d(\G) = {\rm det}\,(-A)$. For a vertex $v$ of $\G$ let
 $\partial_v \G$ denotes the graph obtained from $\G$ by deleting $v$ together
 with all its incident links. If $X$ is a subgraph of $\G$ such that $v \notin
 X$, then $\p_vX$ denotes the subgraph of $\G$ which is obtained from $X$ by
 deleting all the vertices in $X$ closest to $v$ together with their links. In
 what follows we suppose that $\G$ is a tree and $X$ is a subtree; in this
 case there is always the unique vertex in $X$ closest to $v$. Let
 $X_1,\dots,X_N$ be all the non--empty subtrees of $\G$, and put $P = \Z [
 X_1,\dots,X_N]$, where we identify $1$ with the empty subtree and regard the
 disjoint union of subtrees as their product. Then the dicriminant $d$ extends
 to a ring homomorphism $$d\,:\,P \to \Z$$ and $\p_v$ generates a ring
 derivation $$\p_v\,:\,P \to P\,.$$ Denote also $d_v(\G) = d(\p_v \G)$ and
 $d_{vv}(\G) = d(\p_v\p_v \G)$. Let $a_v$ be the weight of a vertex $v \in \G$
. \\

\noin {\bf Proposition 1.}  {\it Let the notation be as above, and let $\G$ be
 a weighted tree. Then

\noin a) For any vertex $v \in \G$ we have
\be d(\G) = -a(v)d_v(\G) - d_{vv} (\G) \,.\ee
b) If $\G$ is a linear tree with the end vertices $v$ and $w$, then
\be d_v(\G)d_w(\G) - d(\G)d_{vw}(\G) = 1\,. \ee
c) Let $[v,\,w]$ be a link of $\G$. Put $\G \setminus \,]\,v,\,w\,[ \,= \G_1
 \cup \G_2$, where $v \in \G_1$ and $w \in \G_2$. Then
\be d(\G) = d(\G_1)d(\G_2) - d_v(\G_1)d_w(\G_2)\,. \ee
d) Let $T$ be a twig of $\G$ incident with a branch vertex $v_0$ of $\G$, and
 let $v$ be the tip of $T$. Put $d_T (\G) = d(\G - T - v_0)$. Then
\be d_T (\G) = d_v(\G)d(T) - d(\G)d_v(T)\,. \ee}
{\bf Corollary.} {\it Let $T$ be a twig of $\G$ such that $d(T) \neq 0$.
 Assume that $d(\G) = 1$. Denote $a = d_T (\G)\, / \, d(T)$. Let $\[a\]$ be
 the least integer bigger than or equal to $a$ and $\]a\[ = \[a\] - a$ be the
 upper fractional part of $a$. Then, in the notation of (d) above, we have}
\be d_v(\G) = \[a\] \,\,\, \,\,\, {\rm and}\,\,\, \,\,\,\ind(T) = \]a\[ \,.
 \ee
\section{Puiseux data as graph discriminants $\,\,\,\,\,\,\,\,\,\,\,\,\,\,\,
\,\,\,$ (after Eisenbud and Neumann)}

Let $(C, \, 0)$ be a germ of an irreducible analytic curve, and let
$$
        x=t^m,\sp\sp y=a_n t^n+a_{n+1}t^{n+1}+...
,\sp\sp a_n\ne 0 \,,
$$
be its analytic parametrization. We may assume that $m<n$ and $m$ does not
 divide $n$. Folowing [A] set $d_1  = m \,,\, m_1  = n;$
%$$ d_i = \gcd( d_{i-1}, m_{i-1}) \,,\sp\sp
%m_i  =  \min \{\, j \mid a_j \ne 0 \,
% \,\,\,\,\,{\rm and}\,\,\,
%\,\,\,j \not\equiv 0 \,
%\,\,({\rm mod}\,d_i)\,\}
%\,,\,\, i>1\,.
%$$
$$ d_i = \gcd( d_{i-1}, m_{i-1}) \,,\sp
 m_i  =  \min \{\, j \mid a_j \ne 0 \,\,\,\,{\rm and}\,\,\,\,
         j \not\equiv 0 \,\,\,({\rm mod}\,d_i)\,\}
\,,\,\, i>1\,.
$$
Let $h$ be such that $d_h \ne 1$, $d_{h+1} =1$.
Thus,  $m_i$ resp. $d_i$ are defined for  $i = 1, ..., h$ resp. for
 $i = 1, ..., h+1$, and
$$
        0< n = m_1<m_2<...<m_h, \sp\sp m=d_1>d_2>...>d_{h+1}=1\,\,.
$$
Set $q_1 = m_1,\,\, q_i = m_i - m_{i-1}$ for $i = 2, ..., h$,
and
\be
        r_i = (q_1 d_1 +...+ q_i d_i)/d_i, \,\,\, i=1,...,h\,\,.  \ee
The sequence $(m; \,m_1,m_2,...,m_h)$ is called {\it the Puiseux
 characteristic sequence of the singularity $(C,\,0)$} [A, Mil]. The whole
 collection $(m_i),\,(d_i),\,(q_i),\,(r_i)$ we call {\it the Puiseux data}. We
 have the following \\

\noin {\bf Proposition 2} [EN]. {\it
a) Let $X \to \C^2$ be the embedded minimal resolution of the singularity
 $(C,\,0)$ with the exceptional divisor $E = \cup E_i$. The proper preimage of
 $C$ in $X$ we denote by the same letter. Then the dual graph $\G_{E,\,C}$ of
 $E$ near $C$ looks like
$$        
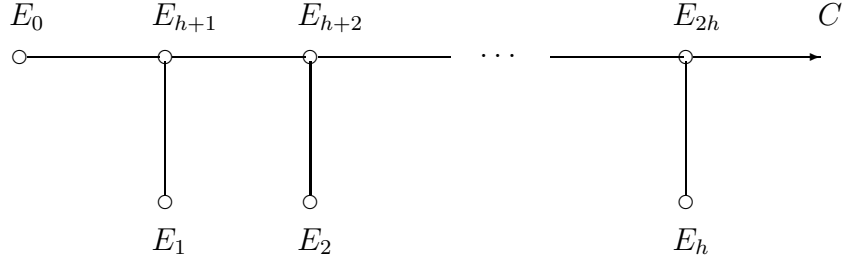
\begin{picture}(1000,90)
          \put(66,82){$E_0$}
          \put(70,70){\circle{5}}
          \put(73,70){\line(1,0){50}}
          \put(120,82){$E_{h+1}$}
          \put(125,70){\circle{5}}
          \put(128,70){\line(1,0){50}}
          \put(175,82){$E_{h+2}$}
          \put(180,70){\circle{5}}
          \put(183,70){\line(1,0){50}}
          \put(244,70){$\ldots$}
          \put(271,70){\line(1,0){50}}
          \put(317,82){$E_{2h}$}
          \put(322,70){\circle{5}}
          \put(325,70){\vector(1,0){48}}
          \put(372,82){$C$}
          %\put(376,70){\circle{5}}
          \put(125,68){\line(0,-1){50}}
          \put(125,15){\circle{5}}
          \put(120,-3){$E_1$}
          \put(180,68){\line(0,-1){50}}
          \put(180,15){\circle{5}}
          \put(175,-3){$E_2$}
          \put(322,68){\line(0,-1){50}}
          \put(322,15){\circle{5}}
          \put(317,-3){$E_h$}
          \end{picture}
  $$
where the edges mean linear chains of vertices of valency two, which are not
 shown. \\

b) Denote by $R_i$, $D_i$ and $S_i$ the connected components of
$\G_{E,\,C}-E_{h+i}$ which are to the left, to the bottom and to the
right of the node $E_{h+i}$, respectively. Denote by $Q_i$ the linear chain
between $E_{h+i-1}$ and $E_{h+i}$ (excluding $E_{h+i-1}$ and $E_{h+i}$).
Then}
$$
        d(R_i)={r_i\over d_{i+1}}\,,\sp
        d(D_i)={d_i\over d_{i+1}}\,,\sp
        d(S_i)=1\,,                 \sp
        d(Q_i)={q_i\over d_{i+1}}\,\, .\sp
$$
This graph is usually called {\it a comb} (see e.g. [FZ1]); M. Miyanishi
 suggested more pleasant name {\it a Christmas tree} (in this case it is drown
 in a slightly different manner).

\section{Expressions of the local BMY-ingredients via the Puiseux data}

\noin {\bf Proposition 3} [OZ, (5.2), (5.4)]. {\it Let $(C,\,0)$ be the local
 branch of $D$ at a cusp $P = P_i$ of $D$. Then in the notation of Sections 6
 and 9 we have
\be \mu = 1 - d_1 + \sum\limits_{i=1}^h r_i ({d_i\over d_{i+1}} -1) = 1 - d_1
 + \sum_{i=1}^h q_i (d_i -1)\,; \ee
\be        2\mu+H_E^2 = -{d_1\over r_1}
        +\sum_{i=1}^h {r_i\over d_{i+1}}({d_i\over d_{i+1}}-{d_{i+1}\over
 d_i}) = -{d_1\over q_1} +\sum_{i=1}^h q_i(d_i-{1\over d_i})\,;\ee
\be       - N_E^2  = \]{d_1\over r_1}\[ + \sum_{i=1}^h \]{r_i\over d_i}\[ =
 \[{d_1\over r_1}\] -  {d_1\over r_1}
        +\sum_{i=1}^h (\[{r_i\over d_i}\] - {r_i\over d_i})\,; \ee
\be        2\mu+(K_E+{\~ D}_E)^2 = -\[{d_1\over r_1}\]
 +\sum_{i=1}^h ({r_i d_i\over d_{i+1}^2}-\[{r_i\over d_i}\])\,\,.
\ee
In particular, if $m$ and $n$ are coprime (i.e.
there is the only one Puiseux characteristic pair), then $\rm (cf.\,\,
 [Mil,\, p. \,95])$ $\mu = (m-1)(n-1)$ and
\be        -H_E^2 = (m-2)(n-2) + (m-n)^2/mn \,,\,\,
\,\,\,\, - N_E^2 = \]{m\over n}\[ + \]{n\over m}\[\,\,\,.\ee}
The proof is based on the formulas in Lemma 3, where the corresponding entries
 $b_{ij}$ have been expressed in terms of the Puiseux data as it is done in
 Proposition 2 above, by using the graph discriminant calculus from Section 8.

\section{Local inequality for irreducible   $\,\,\,\,\,\,\,\,\,\,
\,\,\,\,\,\,\,\,\,\,
\,\,\,\,\,\,\,\,\,\,\,
\,\,\,\,\,\,\,\,$ singularities}

Now we can prove the local inequality (14) for a cusp $P=P_i \in D$. As above,
 we regard the local branch of $D$ at $P=P_i$ as a germ $(C,\,0)$ of an
 analytic curve. \\

\noin {\bf Theorem 2} [OZ, (6.2)]. {\it In the notation as above, for an
 irreducible plane curve germ $(C,\,0)$ one has
\be -(K_E + {\~ D}_E)^2 \le (1 - {1 \over m})\,\mu\,,\ee
where $m$ is the multiplicity and $\mu$ is the Milnor number of $(C,\,0)$.
 The equality in (31) holds iff $m=2$.} \\

The proof proceeds as follows. (31) is equivalent to the inequality
 $$\,\,\,\,\,\,\,\,\,\,\,\,\,\,\,\,\,\,\,\,\,\,\,\,\,\,\,\,\,\,\,\,\,\,
\,\,\,\,\,\,\,\,\,\,\,\,\,\,\,\,\,\,\,\,\,\,\,\,\,\,\,\,\, \mu + (K_E+
 {\~ D}_E)^2 - {\mu\over m} \ge 0\,. \,\,\,\,\,\,\,\,\,\,\,\,\,\,\,\,\,
\,\,\,\,\,\,\,\,\,\,\,\,\,\,\,\,\,\,\,\,\,\,\,\,\,\,\,(31')$$ Using (25)
 and (26--29) we can express the quantity at the left hand side as
\be        \mu + (K_E+ {\~ D}_E)^2 - {\mu\over m} =
        d_1(1-{1\over q_1})-{1\over d_1} + N_E^2
        + \sum_{j=1}^h q_j(1-{d_j\over d_1})(1-{1\over d_j}) \,.\ee
It is easily varified that this quantity vanishes when $m = 2$. The last
 sum in (32) is always positive. Let us show that for $m > 2$ the rest at
 the right hand side of (32) is also positive. Indeed, by (28) we have
$${d_1 \over q_1} - N_E ^2 =  \[{d_1\over q_1}\]
        +\sum_{i=1}^h \]{r_i\over d_i}\[ < \[{m \over n}\] + h\ = 1+h\,.$$
 Thus, it is enough to show that $d_1 - {1 \over d_1} - (1 +h) = m - {1 \over
 m} - (1+h) > 0$. It is true for $m \ge 4$ because $h\le\log_2m$; it is also
 true for $m = 3$ because then $h = 1$, and we are done. \\

\noin {\bf Remark.} The estimate in Theorem 2 is asymptotically sharp in the
 following sense. For any positive integer $m$ and for any $\epsilon>0$ there
 exists an irreducible curve germ $(C,\,0)$ of
multiplicity $m$ such that
$$\,\,\,\,\,\,\,\,\,\,\,\,\,\,\,\,\,\,\,\,\,\,\,\,\,\,\,\,\,\,\,\,\,\,\,\,\,
\,\,\,\,\,\mu + (K_E+ {\~ D}_E)^2 < \mu + H_E^2<(1+\epsilon)\mu/m\,\,. \,\,
\,\,\,\,\,\,\,\,\,\,\,\,\,\,\,\,\,\,\,\,\,\,\,\,\,\,\,\,$$
Indeed, consider the curve {$x^m = y^n$}, where gcd$(m,\,n) = 1$ and $n \gg m$.

\section{Local inequality for arbitrary  singularities}
Here we prove (31) in general case when $(C,\,0)$ is not supposed to be
 irreducible, combining our approach with those of F. Sakai [Sa] (see the
 discussion at the end of Section 3). Let $r$ be the number of local branches
 of $C$, and let $$(m_1=m, \,m_2, \dots, m_n, {\underbrace{1,\dots,1}_{r}})$$
be {\it the multiplicity
 sequence} of $(C,\,0)$, i.e. the sequence of multiplicities of $(C,\,0)$ in
 all its infinitely near points (where $n$ is the total number of blow ups in
the resolution process). Recall [Mil] that \be\mu  + r - 1 =
 \sum\limits_{j=1}^n m_j(m_j - 1)\,.\ee  Remind also that the blow-up at an
 infinitely near point which belongs to only one irreducible component of the
 exceptional divisor is called {\it sprouting} or {\it outer} blow-up, and the
 other blow-ups are called {\it subdivisional} or {\it inner} [MaSa, FZ1].
 Following [Sa] denote by $\omega$ the
 number of subdivisional blow--ups, and set $$ \eta = \sum\limits_{j=1}^n (m_j
 -1)\,.$$
\noin {\bf Lemma 4.} {\it In the notation as above, the following identities
 hold:}
\be -E^2 = \omega\,,\,\,\,\,\,\,EK_E = \omega - 2\,,\,\,\,\,\,\, E^2 + EK_E =
 -2 \,,\,\,\,\,\,\,\, -K_E^2=n \ee
\be ED_E' = r \,,\,\,\,\,\,\, K_ED_E' = \sum\limits_{j=1}^n m_j \,,\,\,\,\,\,
\, K_E(K_E + D_E') = \eta \ee
\be -D_E'^2 = \sum\limits_{j=1}^n m_j^2 \,,\,\,\,\,\,\, -D_E'(K_E + D_E') =
 \sum\limits_{j=1}^n m_j(m_j -1) = 2 \delta\ee
\be \mu + (K_E + \D_E)^2 = (\eta - 1) + (\omega - 1) + (r - 1) \ee
\noin {\bf Proof.} The third equality in (34) resp. the second one in (36)
immediately follows
 from the preceding ones. The first equality in (35) is evident. The other
 formulas in (34) - (36) are proven by an easy induction by the number of
 steps in the resolution process (cf. [MaSa, Lemma 2] for the first equality
 in (34)). To prove (37), transform its left hand side by using (19') and
 (34)--(36) as follows: $$ \mu + (K_E + \D_E)^2 = 1 + (K_E + \D_E)(K_E + E) =
 1 + (K_E + D_E' + E)(K_E + E)$$ $$ = 1 + K_E (K_E + D_E') + 2 EK_E + E^2 +
 ED_E' = \eta + \omega + r - 3\,.$$ \qed \\

The next theorem is a generalization of Theorem 2 to the case when $(C,\,0)$
 is not necessarily irreducible. \\

\noin {\bf Theorem 3.} {\it The inequality (31) is valid for any singular
 plane curve germ $(C,\,0)$, with the equality sign only for an irreducible
 singularity of multiplicity two.} \\

\noin {\bf Proof.} Replace (31) by the equivalent inequality ($31'$). Applying
 (37) we obtain one more equivalent form of (31): \be \eta + \omega + r - 3
 \ge {\mu \over m}\,. \ee This inequality was proven
 in [Sa]. For the sake of completeness we remind here the proof. From (33) it
 follows that \be \eta \ge {\mu + r - 1 \over m}\,,\,\,\,\,\,\,{\rm or}\,\,\,
\,\,\,\eta - {\mu \over m} \ge {r-1 \over m}\,. \ee Therefore, it is enough to
 proof the inequality \be \omega + r - 3 + {r-1 \over m} \ge 0\,, \ee which is
 in turn a consequence of the following one  \be \omega + r \ge 3\,. \ee
 Notice, following [Sa], that $\omega \ge 2$ as soon as at least one
 irreducible branch of $C$ at $0$ is singular, and $\omega = 1$ otherwise.
 But in the latter case $r \ge 2$, because $C$ is assumed being singular.
 This proves (41), and thus also (31). Due to (41) the inequality (40), and
 hence also (31), is strict if $r > 1$. In the case when $r = 1$ by Theorem
 2 the equality sign in (31) corresponds to $m = 2$. This completes the proof.
 \qed

\section{On the rigidity problem for rational $\,\,\,\,\,\,\,\,\,\,\,\,$
 cuspidal plane curves}

Let $Y$ be a a smooth affine algebraic surface $/\C$. Assume that $Y$ is
 $\Q$--acyclic, i.e. $H_i (Y;\,\Q) = 0$ for all $i > 0$, and that $Y$ is of
 log--general type, i.e. ${\bar k}(Y) = 2$. In [FZ1] the problem was posed
 whether such a surface should be rigid. The latter means that $h^1
 (\Theta_X\langle\,\~ D\,\rangle)=0$, where $X$ is a minimal smooth completion
 of $Y$ by a simple normal crossing divisor $\D$ and
 $\Theta_X\langle\,\~ D\,\rangle$ is the logarithmic tangent bundle of $X$
 along $\~ D$. The rigidity holds in all known examples of $\Q$--acyclic
  surfaces of log--general type [FZ1]. Moreover, in all those examples $Y$
  (or, more precisely, the logarithmic deformations of $Y$, see [FZ1]) is
 unobstructed, i.e. $h^2 (\Theta_X\langle\,\~ D\,\rangle)=0$, and therefore
 also the holomorphic Euler characteristic
 $\chi(\Theta_X\langle\,\~ D\,\rangle)$ vanishes (indeed, by Iitaka's Theorem
 [Ii, Theorem 6] $h^0 (\Theta_X\langle\,\~ D\,\rangle)=0$ as soon as $Y$ is of
 log--general type). We have the identity
 $\chi(\Theta_X\langle\,\~ D\,\rangle) = K(K + \D)$ [FZ1, Lemma 1.3(5)], where
 $K = K_X$. Since $\D$ is a curve of arithmetic genus zero [FZ1, Lemma 1.2],
 the equality $\chi(\Theta_X\langle\,\~ D\,\rangle) = K(K + D) = 0$ is
 equivalent to the  following one \be(K + \D)^2 = -2\,.\ee Thus, if $Y$ is
 unobstructed, then it is rigid iff (42) holds.

Consider now an irreducible plane curve $D$. It is easily seen (cf. [Ra]) that
 $Y = \P^2 \setminus D$ is a $\Q$--acyclic surface iff $D$ is a rational
 cuspidal curve. Furthermore, if $D$ has at least three cusps, then $Y$ is of
 log--general type [Wak]. The rigidity of $Y$ is equivalent to $D$ being
 projectively rigid in the following sense: any small deformation of $D$,
 which is a plane rational cuspidal curve with the same types of cusps (i.e.
 an equisingular embedded deformation), is projectively equivalent to $D$
 [FZ2, (2.1)]. Once again, the rigidity holds in all known examples [FZ2,
 (3.3)], as well as the equality in (42) [FZ2, (2.1)].  Here we prove the
 following \\

\noin {\bf Proposition 4.} {\it A projectively rigid rational cuspidal plane
 curve cannot have  more than 9 cusps.} \\

Before giving the proof we remind the notation. Let $\sigma\,:\,X \to \P^2$ be
 the minimal embedded resolution  of singularities of $D$, $K = K_X$ be the
 canonical divisor, $\D = \sigma^{-1}(D)$ and $K+{\~ D} = H+N$ be the Zariski
 decomposition. For a fixed cusp $P \in  {\rm Sing}\,D$ let $K_E + \D_E =
 H_E + N_E$ be the local Zariski--Fujita decomposition, where $E =
 \sigma^{-1}(P)$ is the exceptional divisor.

The proof of Proposition 4 is based on the following two observations.\\

\noin {\bf Lemma 5.} {\it For the negative part $N_E$ of the local
 Zariski--Fujita decomposition over a cusp $P \in D$ the inequality $-N_E^2 >
  1/2$ holds. } \\

\noin {\bf Proof.} From (28) it follows that \footnote{ Recall that $\]a\[$
 denotes $\[a\]-a$, where $\[a\] :=\min\{n\in\Z\,|\,n\ge a\}$.}
\be       - N_E^2  \,= \,\]{d_1\over r_1}\[ +
 \sum_{i=1}^h \]{r_i\over d_i}\[\, \,\le \,\,\]{d_1\over r_1}\[ +
 \]{r_1\over d_1}\[ \,= \, \]{m\over n}\[ + \]{n \over m}\[ \,.\ee By the
 definition of the Puiseux sequence (see section 9) we have $0 < {m \over n}
 < 1$ and ${m \over n} \neq {1 \over 2}$. Thus, the desired inequality follows
 from (43) and the next estimate, which is an easy exercise. \qed \\

\noin {\bf Claim.} {\it If $0< x <1$, then $\]x\[ + \]{1 \over x}\[ \ge {1
 \over 2}$, where the equality holds only for $x = {1 \over 2}$. }\\

\noin {\bf Lemma 6.} {\it Let $D$ be a rational cuspidal plane curve with at
 least three cusps. Then in the notation as above we have \be H = (d-3)L' +
 \sum_{P \in {\rm Sing}\,D} H_E\,\,\,{\rm and }\,\,\,N = \sum_{P \in {\rm
 Sing}\,D} N_E,\,\ee i.e. the global Zariski decomposition agrees with the
 local Zariski--Fujita ones.} \\

\noin {\bf Proof.} By [Wak] we have ${\bar k}(Y) = 2$, where $Y = \P^2
 \setminus D = X \setminus \D$. Thus, being a smooth $\Q$--acyclic surface of
 log--general type, $Y$ does not contain any simply connected curve (this was
 first proven in [Za] for acyclic surfaces and then generalized in [MT] to
 $\Q$--acyclic ones). In particular, $X$ does not contain any $(-1)$--curve
 $C$ with $C\cdot \D = 1$. Since $D$ has at least three cusps, the dual graph
 $\G_{\D}$ of $\D$ has at least three branching points. Under these conditions
 the lemma follows from the results in [Fu, (6.20-6.24)] (see also [OZ,
 Theorem 1.2]).  \qed \\

\noin {\bf Proof of Proposition 4.} Let $\kappa$ be the number of cusps of
 $D$. Evidently, we may suppose that $\kappa \ge 3$. It follows from Lemmas 5
 and 6 that
$$(K+{\~ D})^2 = H^2 + N^2 = H^2 + \sum_{P \in {\rm Sing}\,D} N_E^2 < H^2
- {1 \over 2}\,\kappa\,.$$
Due to BMY-inequality (13) we also have $H^2 \le 3$, and hence
$$(K+{\~ D})^2 < 3 - {1 \over 2}\,\kappa\,.$$
Set $h^i = h^i (\Theta_X\langle\,\~ D\,\rangle)\,,\,\,i = 0,\,1,\,2$. The
 surface $Y =  \P^2 \setminus D$ being of log--general type [Wak], by Iitaka's
 Theorem [Ii, Theorem 6] we have $h^0 = 0$. Since $D$ is assumed to be rigid,
 i.e. $h^1 = 0$, we also have $\chi(\Theta_X\langle\,\~ D\,\rangle) = h^2 =
 K(K+{\~ D}) \ge 0$, i.e. $(K+{\~ D})^2 \ge -2$. It follows that
\be \kappa < 6 - 2(K+{\~ D})^2 \le 10 \,,\ee which completes the proof. \qed
 \\

\noin {\bf Remark.} Actually, for a rational cuspidal plane curve with at
 least three cusps we have proved the inequality
\be \kappa < 6 - 2(K+{\~ D})^2 = 10 - 2K(K+{\~ D}) \,.\ee
Therefore, $\kappa < 10$ as soon as $K(K+{\~ D}) \ge 0$, which is the case if
 $D$ is rigid.

\vs

\begin{center} {\bf REFERENCES}
\end{center}

\vs

{\footnotesize

\noin [A]  S. S. Abhyankar. Expansion  technique in algebraic geometry.
 {\it Tata Inst. of Fund. Res.}, Bombay, 1977 \\

\noin [BK] E. Briskorn, H. Kn\"orrer. Plane algebraic curves. {\it
 Birkh\"auser-Verlag,} Basel e.a., 1986 \\

\noin [DG] D. Drucker, D.M. Goldschmidt. Graphical evaluation of
sparce determinants. {\it Proc. Amer. Math. Soc.}, {\bf 77}  (1979), 35 - 39\\

\noin [EN] D. Eisenbud, W. D. Neumann. Three-dimensional link theory and
 invariants of plane curve singularities. {\it Ann.Math.Stud.}  {\bf  110},
 {\it Princeton Univ. Press}, Princeton 1985\\

\noin [FZ1] H. Flenner, M. Zaidenberg. $\bf Q$--acyclic surfaces and their
 deformations. {\it Contemporary Mathem.} {\bf 162} (1964), 143--208 \\

\noin [FZ2] H. Flenner, M. Zaidenberg. On a class of rational cuspidal plane
 curves. {\it Preprint} (1995), 1--28 \\

\noin [Fu] T. Fujita. On the topology of non-complete algebraic surfaces, {\it
    J. Fac. Sci. Univ. Tokyo (Ser 1A)}, {\bf 29} (1982), 503--566 \\

\noin [Hi] A. Hirano. Construction of plane curves with cusps, {\it Saitama
 Math. J.} {\bf 10} (1992), 21--24 \\

\noin [H] F. Hirzebruch. Singularities of algebraic surfaces and
 characteristic numbers, {\it The Lefschetz Centennial Conf. Part.I
 (Mexico City 1984), Contemp. Math.} {\bf  58} (1985),  141-155 \\

\noin [Ii] Sh. Iitaka. On logarithmic Kodaira dimension of algebraic
 varieties. In: {\it Complex Analysis and Algebraic Geometry, Cambridge Univ.
 Press}, Cambridge e.a., 1977, 175--190\\

\noin [Iv] K. Ivinskis, Normale Fl\"achen und die Miyaoka--Kobayashi
 Ungleichung. {\it Diplomarbeit}, Bonn, 1985 \\

\noin [KoNaSa] R. Kobayashi, S. Nakamura, F. Sakai. A numerical
 characterization of ball quotients for normal surfaces with branch loci.
 {\it Proc. Japan Acad.} {\bf  65(A)} (1989), 238--241 \\

\noin [Ko] R. Kobayashi. An application of K\"ahler--Einstein metrics to
 singularities of plane curves. {\it Advanced Studies in Pure Mathem., Recent
 Topics in Differential and Anal. Geom.} {\bf 18-I} (1990), 321--326 \\

\noin [MaSa] T. Matsuoka, F. Sakai. The degree of rational cuspidal curves.
 {\it Math. Ann.} {\bf 285} (1989), 233--247\\

\noin [Mil] J. Milnor. Singular points of complex hypersurfaces. {\it
 Ann.Math.Stud.} {\bf  61}, {\it Princeton Univ. Press},  Princeton, 1968 \\

\noin [MT] M. Miyanishi, S. Tsunoda. Abscence of the affine lines on the
 homology planes of general type. {\it J. Math. Kyoto Univ.} {\bf 32} (1992),
 443--450 \\

\noin [Miy] Y. Miyaoka. The minimal number of quotient singularities on
 surfaces with given numerical invariants, {\it Math. Ann.} {\bf  268}
 (1984),  159--171 \\

\noin [Na] M. Namba. Geometry of projective algebraic curves. {\it Marcel
 Dekker},  N.Y. a.e., 1984 \\

\noin [Ne] W.D. Neumann. On bilinear forms represented by trees. {\it Bull.
 Austral. Math. Soc.} {\bf  40}  (1989),  303-321 \\

\noin [OZ] S.Y. Orevkov, M.G. Zaidenberg. Some estimates for plane cuspidal
 curves. In: {\it Journ\'ees singuli\`eres et jacobiennes, Grenoble 26--28 mai
 1993.}  Grenoble, 1994, 93--116 (see also Preprint MPI/92-63, 1992, 1--13)\\

\noin [Ra] C.P. Ramanujam. A topological characterization of the affine
plane as an algebraic variety. {\it Ann. Math.} 94 (1971), 69-88 \\

\noin [Sa] F. Sakai. Singularities of plane curves. {\it Preprint} (1990),
 1-10\\

\noin [Va] A.N. Varchenko. Asymptotics of integrals and Hodge structures. In:
 {\it Itogi Nauki i Techniki, Series "Contempor. Problems in Mathem."}
 {\bf 22} (1983), 130--166 (in Russian)\\

\noin [Wak] I. Wakabayashi. On the logarithmic Kodaira dimension of the
 complement of a curve in $\P^2$. {\it Proc. Japan Acad.} {\bf 54(A)} (1978),
  157--162 \\

\noin [Wal] R. J. Walker.  Algebraic curves. {\it Princeton Univ. Press},
  Princeton, 1950 \\

\noin [Yo1] H. Yoshihara. Plane curves whose singular points are cusps.
 {\it Proc. Amer. Math. Soc.} {\bf 103} (1988), 737--740 \\

\noin [Yo2] H. Yoshihara. Plane curves whose singular points are cusps and
 triple coverings of $\P^2$.  {\it Manuscr. Math.} {\bf  64} (1989), 169-187
 \\

\noin [Za] M. Zaidenberg. Isotrivial families of curves on affine surfaces and
 characterization of the affine plane. {\it Math. USSR Izvestiya} {\bf 30}
 (1988), 503-531; Addendum, {\it ibid.} {\bf 38} (1992), 435--437

\vs

\noin Stepan Orevkov\\

System Research Institute RAN\\
Moscow, Avtozavodskaja 23, Russia\\
e-mail: OREVKOV@GLAS.APS.ORG
\vs

\noin Mikhail Zaidenberg\\

Universit\'{e} Grenoble I \\
Laboratoire de Math\'ematiques associ\'e au CNRS\\
BP 74\\
38402 St. Martin d'H\`{e}res--c\'edex, France\\
e-mail: ZAIDENBE@FOURIER.GRENET.FR}

\end{document}